\documentclass[12pt,a4paper]{article}
 \usepackage{emlines2,amsmath,amstext,amssymb,amscd}

\begin{document}
\title{The octahedron recurrence and RSK-correspondence}
\author{V. I. Danilov and G.A. Koshevoy\thanks{Central Institute
of Economics and Mathematics, Russian Academy of Sciences,
Nakhimovskij Prospekt 47, Moscow 117418, Russia. Email:
danilov@cemi.rssi.ru, koshevoy@cemi.rssi.ru. This research was
supported in part by grants NSch-6417.2006.6, NWO and RFBR
047.011.204.017, and CNRS and RFBR CNRS\,a 05-01-02805.}}
\date{14.03.07}
\maketitle

\section{Introduction}

In \cite{umn} we proposed a new approach to the RSK
(Robinson-Schensted-Knuth) correspondence based on bi-crystal
structure of arrays. However, it was not easy to compare our
construction with the classical one. The point is that the
classical RSK-correspondence uses the bumping (or sliding)
procedure (see, for example \cite{Fulton,Stan-2}), while our array
construction is done via condensation operations. Recent papers
\cite{BZ,kirill,NY,oct} show  connections of this issue with the
Dodgson rule for calculation of determinants. Roughly speaking,
the RSK-correspondence is a tropicalization of the Dodgson rule.
To make this statement rigorous is one of the goals of this paper.

We start with an ``algebraic'' RSK-correspondence due to Noumi and
Yamada \cite{NY}. Given a matrix $X$, we consider a pyramidal
array of solid minors of $X$. It turns out that this array
satisfies an algebraic variant of octahedron recurrence. The main
observation is that this array can also  be constructed with the
help of some square `genetic' array. For example, if a `genetic'
array is positive then the corresponding matrix $X$ is totally
positive. Furthermore, any totally positive matrix can be obtained
by this way (\cite{BZ}).

Next we tropicalize this algebraic construction and consider
$T$-{\em polarized} pyramidal arrays (that is arrays satisfying
octahedral relations). As a result we get several bijections, viz:
a) a linear bijection between non-negative arrays and supermodular
functions; b) a piecewise linear bijection between supermodular
functions and the so called infra-modular functions; c) a linear
bijection between infra-modular functions and plane partitions. A
composition of these bijections yields a bijection between
non-negative arrays and plane partitions coinciding with the
modified RSK-correspondence defined in \cite{umn}.

The paper is organized as follows. In Section 2 we present an
algebraic version of the octahedron recurrence (OR). In Section 3
we consider tropically polarized functions and their relation to
supermodular functions. In Section 4, a construction of tropically
polarized functions via the tropical `genetic' array is given. In
Section 5, we show (Theorem 3) that the tropical OR gives a
bijection between the set of supermodular functions and that of
infra-modular functions. Section 6 contains a proof of Theorem 3.

We would like to thank the referee for attracting our attention to
the paper \cite{NY}, and proposing to consider our `octahedron'
RSK-correspondence (\cite{arxiv}) as a tropicalization of the
`algebraic' RSK-correspondence due to Noumi and Yamada. We thank
F. Zak for his numerous remarks and improvements.

\section{Algebraic case}

Although we are interested in two-dimensional arrays of numbers,
relations between them will be constructed via three-dimensional
arrays. Therefore, we shall operate also with functions on
three-dimensional lattices.

Specifically, consider an $n\times n$-matrix $X=(X(r,c))$ with
rows $r=1,...,n$ and columns $c=1,...,n$. (A general case of
rectangular matrices can be reduced to this case without major
difficulties. However, for simplicity, we confine ourselves to the
case of square matrices.) We are interested in all solid minors of
$X$. We arrange these minors in a three-dimensional pyramidal
array.

Let $ABCDE$ be a pyramid with base  $ABCD$, whose vertices
$A=(0,0,0)$, $B=(2n,0,0)$, $C=(0,2n,0)$, $D=(2n,2n,0)$ are located
at `level' $0$, and the top vertex $E=(n,n,n)$.

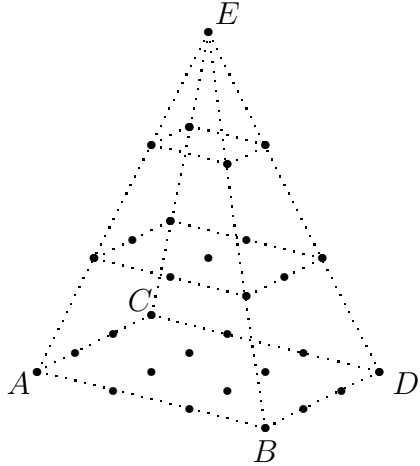
\unitlength=.5mm \special{em:linewidth 0.4pt}
\linethickness{0.4pt}
\begin{picture}(109.00,135.00)(-60,-10)
\put(15.00,20.00){\circle*{2.00}} \put(75.00,5.00){\circle*{2.00}}
\put(105.00,20.00){\circle*{2.00}}
\put(45.00,35.00){\circle*{2.00}}
\put(30.00,50.00){\circle*{2.00}}
\put(70.00,40.00){\circle*{2.00}}
\put(90.00,50.00){\circle*{2.00}}
\put(50.00,60.00){\circle*{2.00}}
\put(60.00,110.00){\circle*{2.00}}
\put(45.00,80.00){\circle*{2.00}}
\put(65.00,75.00){\circle*{2.00}}
\put(75.00,80.00){\circle*{2.00}}
\put(55.00,85.00){\circle*{2.00}}
\bezier{40}(15.00,20.00)(37.00,65.00)(60.00,110.00)
\bezier{24}(15.00,20.00)(45.00,13.00)(75.00,5.00)
\bezier{13}(75.00,5.00)(90.00,13.00)(105.00,20.00)
\bezier{42}(75.00,5.00)(67.00,55.00)(60.00,110.00)
\bezier{40}(60.00,110.00)(82.00,65.00)(105.00,20.00)
\bezier{24}(45.00,35.00)(72.00,28.00)(105.00,20.00)
\bezier{13}(45.00,35.00)(30.00,28.00)(15.00,20.00)
\put(30.00,50.00){\circle*{2.00}}
\put(30.00,50.00){\circle*{2.00}}
\put(50.00,45.00){\circle*{2.00}}
\put(70.00,40.00){\circle*{2.00}}
\put(40.00,55.00){\circle*{2.00}}
\put(50.00,60.00){\circle*{2.00}}
\put(45.00,80.00){\circle*{2.00}}
\put(55.00,85.00){\circle*{2.00}}
\put(65.00,75.00){\circle*{2.00}}
\put(60.00,50.00){\circle*{2.00}}
\put(70.00,55.00){\circle*{2.00}}
\put(80.00,45.00){\circle*{2.00}}
\put(90.00,50.00){\circle*{2.00}}
\put(75.00,80.00){\circle*{2.00}}
\put(60.00,110.00){\circle*{2.00}}
\put(15.00,20.00){\circle*{2.00}}
\put(25.00,25.00){\circle*{2.00}}
\put(35.00,30.00){\circle*{2.00}}
\put(45.00,35.00){\circle*{2.00}}
\put(35.00,15.00){\circle*{2.00}}
\put(45.00,20.00){\circle*{2.00}}
\put(55.00,25.00){\circle*{2.00}}
\put(65.00,30.00){\circle*{2.00}}
\put(55.00,10.00){\circle*{2.00}}
\put(65.00,15.00){\circle*{2.00}}
\put(75.00,20.00){\circle*{2.00}}
\put(85.00,25.00){\circle*{2.00}} \put(75.00,5.00){\circle*{2.00}}
\put(85.00,10.00){\circle*{2.00}}
\put(95.00,15.00){\circle*{2.00}}
\put(105.00,20.00){\circle*{2.00}}
\put(10.00,17.00){\makebox(0,0)[cc]{$A$}}
\put(75.00,-1.00){\makebox(0,0)[cc]{$B$}}
\put(112.00,17.00){\makebox(0,0)[cc]{$D$}}
\put(42.00,39.00){\makebox(0,0)[cc]{$C$}}
\put(65.00,115.00){\makebox(0,0)[cc]{$E$}}
\bezier{30}(60.00,110.00)(53.00,72.00)(45.00,35.00)
\bezier{16}(30.00,50.00)(50.00,45.00)(70.00,40.00)
\bezier{8}(70.00,40.00)(80.00,45.00)(90.00,50.00)
\bezier{16}(90.00,50.00)(70.00,55.00)(50.00,60.00)
\put(50.00,60.00){\circle*{2.00}}
\put(55.00,85.00){\circle*{2.00}}
\bezier{8}(30.00,50.00)(40.00,55.00)(50.00,60.00)
\bezier{8}(45.00,80.00)(55.00,78.00)(65.00,75.00)
\bezier{4}(65.00,75.00)(70.00,77.00)(75.00,80.00)
\bezier{8}(75.00,80.00)(65.00,83.00)(55.00,85.00)
\bezier{4}(55.00,85.00)(50.00,83.00)(45.00,80.00)
\end{picture}

\noindent Denote by $\Pi$ the subset of the pyramid consisting of
integer points $(i,j,k)$, where $i$ and $j$ are equal to $k$
modulo $2$.

Given a matrix $X$ we define a function $F=F_X$ on $\Pi$. On level
zero, $F$ is equal to 1. Elements of $X$ are placed in the points
on level one. More precisely, the number $X(r,c)$ is  assigned to
the point $(2r-1,2c-1,1)$. At level two we place the minors of $X$
of size $2$ and so on. Let $(i,j,k)\in \Pi$, so that $k\le i,j\le
2n-k$. Consider the sub-matrix of $X$ formed by $k$ consecutive
rows $(i-k)/2+1, (i-k)/2+2,..., (i-k)/2+k$ and $k$ consecutive
columns $(j-k)/2+1, (j-k)/2+2,..., (j-k)/2+k$. We set $F(i,j,k)$
to be equal to the determinant of this sub-matrix. To visualize
the construction, one has to consider a pyramidal cone (congruent
to $ABCDE$) with apex at the point $(i,j,k)\in\Pi$. This cone
meets level one along a sub-matrix of $X$. The determinant of this
sub-matrix is assigned to the point $(i,j,k)$. For example,
$F(n,n,n)=\det(X)$.

It is useful to consider the restriction of the function $F$ to
the  faces of the pyramid. Suppose that a point $(i,j,k)\in \Pi$
lies on the face $ABE$, so that $i\ge j=k$. Then the corresponding
minor is formed by the first  $k$ rows and the columns $(j-k)/2+1,
(j-k)/2+2,..., (j-k)/2+k$. Similarly, at the face $ACE$, the
values of $F$ consist of solid  minors of sub-matrices whose
consecutive columns start from $1$.

If we put together values of $F$ on the faces  $ABE$ and $ACE$, we
obtain an $(n+1)\times (n+1)$ array. More precisely, consider the
map
$$
\alpha:\{0,1,...,n\}\times \{0,1,...,n\} \to \Pi,
$$
$$
            \alpha (i,j)=(2i-\min(i,j),2j-\min(i,j),\min(i,j)).
           $$
Put $G=\alpha^*(F)$. For  $j\le i$,  $G(i,j)$ is equal to the
minor of  $X$ formed by the first $j$ rows and columns
$i-j+1,\ldots ,i$. For  $i\ge j$,  $G(i,j)$ is equal to the minor
with the rows $j-i+1,\ldots ,j$ and the first $i$ columns. Note,
that $G(i,j)=1$ if $ij=0$, so actually we get an $n\times n$
array. (We denote this function by $G$ in honor of Gauss since its
construction is related to the Gauss method of elimination of
unknowns.)

 For a generic matrix $X$, the array $G$
determines $X$. (By generic matrix we mean here a matrix with
non-zero solid minors.)

Analogously, one can map the grid $\{0,1,...,n\}\times
\{0,1,...,n\}$ to a union of faces  $DBE$ and  $DCE$ according to
the rule
$$
\beta (\tilde i,\tilde j)=(2n-2\tilde i+\min( \tilde i,\tilde
j),2n-2\tilde j+\min(\tilde i,\tilde j),\min(\tilde i,\tilde j)).
$$
The function $H=\beta^*(F)$  corresponds to the array consisting
of solid minors of  $X$ adjacent to the south-east boundary of the
matrix. For generic matrix $X$, $H$ also determines $X$.

In particular, the bijection between the $G$-data and $H$-data can
be considered as a (birational) mapping from $\mathbb{R}^{n\times
n}$ to $\mathbb{R}^{n\times n}$.

A more explicit transformation of  $G$  into  $H$ (and vice versa)
is based on the following property of the function $F$ called the
Dodgson condensation rule (see, for example \cite{Sp}). Namely,
for $(i,j,k)\in \Pi$ and $k\ge 1$, there holds
\begin{equation}\label{dodg}
F(i,j,k+1)F(i,j,k-1)=F(i-1,j-1,k)F(i+1,j+1,k)-F(i-1,j+1,k)F(i+1,j-
1,k).
\end{equation}
(The values of the function can be taken in any field or
ring.)\medskip

{\bf Definition}. A function   $F:\Pi\to\mathbb R$ is said to be
{\em algebraic polarized} (or $A$-polarized) if $F(*,*,0)=1$ and
the relations (\ref{dodg}) are satisfied.\medskip

Thus, for any $A$-polarized function, its values at the vertices
of each `elementary' octahedron centered at $(i,j,k)$, $i,j =k+1
\mod 2$, satisfy the relation (\ref{dodg}). By the reason, one
often speak about the octahedron recurrence. If a function $F$ is
defined at any five vertices of an elementary octahedron, then its
sixth value can be determined by (\ref{dodg}). (Of course, for
generic case.) This gives a hint that an $A$-polarized function is
determined by its values at some smaller subsets of $\Pi$. For
example, an $A$-polarized function is determined by its values at
the first level, $F(2r-1,2s-1,1)$, $r,s=1,\ldots , n$ (that is by
the matrix $X$). As well, $F$ is determined by its values at any
pair of adjoint faces (for example at faces $ABE$ and $ACE$, where
it is defined by the function $G$).\medskip

Now we discuss another construction of $A$-polarized functions via
`genetic' arrays.

Let us fix an array  $W=(W(i,j), \ 1\le i,j\le n)$. Pick a pair
$(i,j)$. A {\em path} from the column $i$ to the row $j$ is a
sequence $\gamma $ of integer pairs $(a_1 ,b_1 ),...,(a_s ,b _s)$,
such that $(a_1 ,b_1 )=(i,1)$, $(a_s ,b_s )=(1,j)$, and for any
$t=2,\ldots, s$, $(a_t, b_t)-(a_{t-1}, b_{t-1})$ equals either
$(0,1)$ or  $(-1,0)$.

\unitlength=.7mm \special{em:linewidth 0.4pt}
\linethickness{0.4pt}
\begin{picture}(97.00,58.00)(-20,0)
\put(45.00,5.00){\circle*{1.00}} \put(55.00,5.00){\circle*{1.00}}
\put(65.00,5.00){\circle*{1.00}} \put(75.00,5.00){\circle*{1.00}}
\put(85.00,5.00){\circle*{1.00}} \put(45.00,15.00){\circle*{1.00}}
\put(55.00,15.00){\circle*{1.00}}
\put(65.00,15.00){\circle*{1.00}}
\put(75.00,15.00){\circle*{1.00}}
\put(85.00,15.00){\circle*{1.00}}
\put(45.00,25.00){\circle*{1.00}}
\put(55.00,25.00){\circle*{1.00}}
\put(65.00,25.00){\circle*{1.00}}
\put(75.00,25.00){\circle*{1.00}}
\put(85.00,25.00){\circle*{1.00}}
\put(45.00,35.00){\circle*{1.00}}
\put(55.00,35.00){\circle*{1.00}}
\put(65.00,35.00){\circle*{1.00}}
\put(75.00,35.00){\circle*{1.00}}
\put(85.00,35.00){\circle*{1.00}}
\put(45.00,45.00){\circle*{1.00}}
\put(55.00,45.00){\circle*{1.00}}
\put(65.00,45.00){\circle*{1.00}}
\put(75.00,45.00){\circle*{1.00}}
\put(85.00,45.00){\circle*{1.00}}
\put(75.00,6.00){\vector(0,1){8.00}}
\put(74.00,15.00){\vector(-1,0){8.00}}
\put(65.00,16.00){\vector(0,1){8.00}}
\put(65.00,26.00){\vector(0,1){8.00}}
\put(64.00,35.00){\vector(-1,0){8.00}}
\put(55.00,36.00){\vector(0,1){8.00}}
\put(54.00,45.00){\vector(-1,0){8.00}}
\put(46.00,5.00){\vector(1,0){46.00}}
\put(45.00,5.00){\vector(0,1){44.00}}
\put(97.00,3.00){\makebox(0,0)[cc]{$i$}}
\put(48.00,52.00){\makebox(0,0)[cc]{$j$}}
\end{picture}

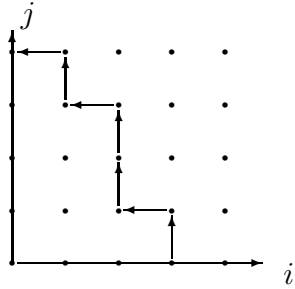
\hfill {\small A path $\gamma$ from the 4-th column to the 5-th
row.} \hfill \bigskip

For a path $\gamma$, we denote by $W^\gamma $ the product $\prod
_{t=1,...,s}W(a_t, b_t )$. Then we set
$$
X(i,j)=F(2i-1,2j-1,1):=\sum _\gamma W^\gamma ,
$$
where $\gamma $ runs over the set of all paths from  $i$ to $j$.
More generally, let us consider  a point $(i,j,k)\in\Pi$. It
determines $k$ consecutive rows
$i_1=(i-k)/2+1,...,i_k=(i-k)/2+k=(i+k)/2$ and  $k$ consecutive
columns  $j_1=(j-k)/2+1,...,j_k=(j+k)/2$. We set
$$
F(i,j,k)=\sum_{(\gamma_1  ,...,\gamma _k )}W^{\gamma _1} \cdots
W^{\gamma _k},
$$
where $(\gamma _1 ,...,\gamma _k)$  runs over the set of $k$
non-intersecting paths from $i_1$ to $j_1 $, from $i_2$ to $j _2$,
$\ldots$, from $i_k$  to $j _k$. As before, $F(i,j,0)=1$.\medskip

{\bf Theorem 1.} {\em  The function $F$ is $A$-polarized.}\medskip

The proof is a combination of the Dodgson rule and the
Lindstr\"om-Gessel-Viennot theorem \cite{Stan-1}, Theorem 2.7.1.
The theorem shows that $F(i,j,k)$ is indeed equal to the minor of
the matrix $X=(X(i,j))$ `illuminated' from the point $(i,j,k)$.
More generally, Lindstr\"om theorem asserts that, for an
increasing tuple of rows $i_1 ,...,i_k$  and an increasing tuple
of columns $j_1 ,...,j_k$, the corresponding minor of the matrix
$X$ is equal to
$$
\sum _{(\gamma_1  ,...,\gamma _k ) }W^{\gamma _1}\cdots W^{\gamma
_k},
$$
where the tuple $(\gamma _1 ,...,\gamma _k)$ runs over the set of
$k$ non-intersecting paths from $i_1$ to $j_1 $, from $i_2$ to $j
_2$, $\ldots$, from $i_k$  to $j _k$.\medskip

Let $F$ be the $A$-polarized function constructed via  a `genetic'
array $W$. We are interested in its values at the faces $ABE$ and
$ACE$. Specifically, we are interested in relations between  the
function $G=\alpha^*(F)$ and the array $W$.

Consider a pair $(i,j)$ with $i\ge j$. Then $G(i,j):=F(2i-j,j,j)$
is equal to the minor of $X$ formed by the columns $i_1
=i-j+1,...,i_j =i$ and rows $1,...,j$. But there is only one
collection of $j$ non-intersecting paths from columns $i_1
=i-j+1,...,i_j =i$ to rows  $1,...,j$, and these paths cover the
whole rectangle $\{1,...,i\}\times \{1,...,j\}$. Therefore, the
corresponding minor of the matrix  $X$ is equal to a product of
$W(a,b)$, $1\le a \le i$, $1\le b \le j$.

From this, we immediately obtain that  $F$ and even $G$ determines
$W$. Namely, we have

\begin{equation}\label{gen1}
            W(i,j)=\frac{G(i,j)G(i-1,j-1)}{G(i-1,j)G(i,j-1)}
\end{equation}
(here we assume that all values $G(*,*)$ are invertible).


For totally positive initial matrices (more precisely, for
matrices with strictly positive minors) we get (after \cite{BZ})
the following diagram of bijections\medskip

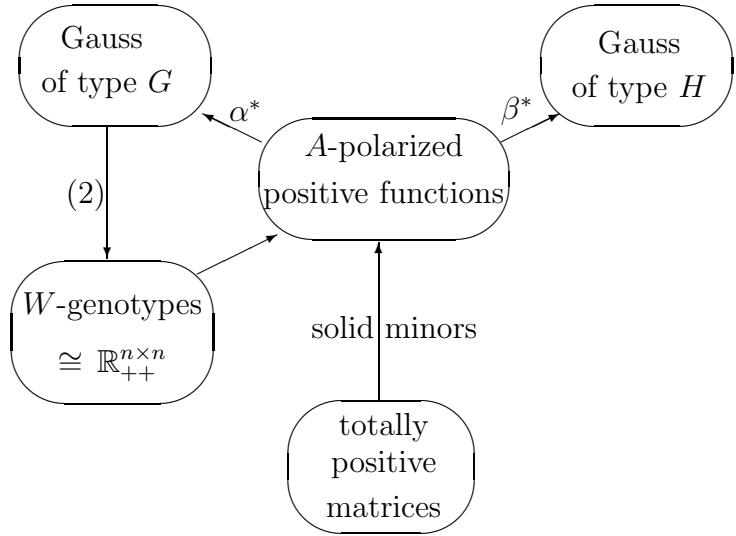
\unitlength=0.7mm \special{em:linewidth 0.4pt}
\linethickness{0.4pt}
\begin{picture}(141.00,115.00)(-20,5)
\put(78,82){\oval(47.00,23.00)[]}
\put(78.00,88.00){\makebox(0,0)[cc]{$A$-polarized}}
\put(78.00,79.00){\makebox(0,0)[cc]{positive functions}}
\put(125.50,103.50){\oval(36.00,23.00)[]}
\put(27.50,103.50){\oval(36.00,23.00)[]}
\put(27,53){\oval(38.00,27.00)[]}
\put(126.00,108.00){\makebox(0,0)[cc]{Gauss}}
\put(126.00,99.00){\makebox(0,0)[cc]{of type $H$}}
\put(25.00,109.00){\makebox(0,0)[cc]{Gauss}}
\put(26.00,100.00){\makebox(0,0)[cc]{of type $G$}}
\put(27.00,58.00){\makebox(0,0)[cc]{$W$-genotypes}}
\put(27.00,47.00){\makebox(0,0)[cc]{$\cong\,\mathbb R^{n\times
n}_{++}$}}

\put(77.50,27.50){\oval(35.00,25.00)[]}
\put(78.00,35.00){\makebox(0,0)[cc]{totally}}
\put(78.00,28.00){\makebox(0,0)[cc]{positive}}
\put(78.00,20.00){\makebox(0,0)[cc]{matrices}}

\put(55.00,89.00){\vector(-2,1){11}}
\put(100.00,89.00){\vector(2,1){11}}
\put(43.00,64.00){\vector(2,1){15}}
\put(77.00,40.00){\vector(0,1){30.00}}
\put(26.00,92.00){\vector(0,-1){25.00}}
\put(52.00,95.00){\makebox(0,0)[cc]{$\alpha^*$}}
\put(103.00,95.00){\makebox(0,0)[cc]{$\beta^*$}}
\put(22.00,79.00){\makebox(0,0)[cc]{$(2)$}}
\put(80.00,54.00){\makebox(0,0)[cc]{solid minors}}
\end{picture}

If we replace totally positive matrices by arbitrary matrices, we
obtain  birational isomorphisms between the corresponding parts of
this diagrams. In \cite{kirill,NY}, the bijection between the set
of $W$-genotypes ($\cong\,\mathbb R^{n\times n}_{++}$) and the set
of Gauss data of type $H$ was considered as an
algebraic\footnote{Such a bijection was called in \cite{NY}
tropical RSK-correspondence, but it seems that using this term in
the algebraic framework is a bit misleading.} RSK-correspondence.

In the next sections we tropicalize  these constructions.

\section{Tropically polarized functions}

We again consider functions on the subset $\Pi$ of the pyramid
$ABCDE$. To distinguish from the algebraic case, we will use small
letters to denote functions in the tropical (or combinatorial)
case.\medskip

{\bf Definition.} A function $f:\Pi\to\mathbb R$ is said to be
 {\em tropically polarized} (shortly, $T$-polarized) if $f(*,*,0)=0$ and
\[
f(i-1,j-1,k)+f(i+1,j+1,k)= \]
\begin{equation}\label{trop}
\max(f(i,j,k-1)+f(i,j,k+1), f(i-1,j+1,k)+f(i+1,j-1,k))
\end{equation}
for any $i$ and $j$ compatible with $k+1$ modulo $2$.\medskip

{\bf  Example 1.} Consider the function $q:\Pi\to\mathbb R$ given
by $q(i,j,k)=k$. Obviously this function is $T$-polarized.\medskip

{\bf Example 2.} Consider the following function $p$ on $\Pi$
$$
                        p(i,j,k)=(i+j-1)k.
$$
It is $T$-polarized. In fact, for an elementary octahedron around
$(i,j,k)$ (where $i-k,j-k$ are odd numbers) the sum of values of
$p$ at end points of any diagonal of this octahedron is equal to
$2(i+j-1)k$.

Note that, for a $T$-polarized function $f$, the function
$f+\alpha p$ is $T$-polarized for any real $\alpha $.\medskip

The next example generalizes the preceding ones.\medskip

{\bf Example 3.} Let $\phi $ and $\psi $ be functions defined on
$\{0,2,...,2n\}$. Consider the following function on $\Pi$
$$
f(i,j,k)=\phi (i-k)-\phi (i+k)+\psi (j-k)-\psi (j+k).
$$

We claim that  $f$ is tropically polarized function. In fact,
$f(i-1,j-1,k)+f(i+1,j+1,k)$ is equal to
\begin{equation}\label{*}
\phi (i-k-1)-\phi (i+k-1)+\psi (j-k-1)-\psi (j+k-1)+ \phi
(i-k+1)-\phi (i+k+1)+\psi (j-k+1)-\psi (j+k+1)
\end{equation}
One can check that $f(i,j,k-1)+f(i,j,k+1)$ is also equal to
(\ref{*}), as well as $f(i-1,j+1,k)+f(i+1,j-1,k)$. Thus, for an
elementary octahedron, the values of $f$ at end points of any of
its diagonal are the same and the relations (\ref{trop}) are
fulfilled. Since $f(i,j,0)=0$, $f$ is $T$-polarized.
$\Box$\medskip

Denote by $Pol$ the set of $T$-polarized functions on $\Pi$. The
set $Pol$ is a polyhedral complex of cones. A cone of the complex
is specified by choosing a cutting of each elementary octahedron
into two half-octahedron such that both halves contain the {\em
propagation vector} $(2,2,0)$. To wit, each elementary octahedron
is determined by its center $(i,j,k)$, where $0<k<n$, $k\le i,j\le
2n-k$ and, modulo 2, $i$ and $j$ differ from $k$. Denote by $\Pi'$
the set of centers of elementary octahedra. The end points of the
diagonal parallel to the propagation vector are the octahedron
vertices $(i-1,j-1,k)$ and $(i+1,j+1,k)\in \Pi$. Through this
diagonal we can make either a vertical cut ($\uparrow$) or a
horizontal cut ($\to$) in order to split the octahedron into
halves. To each function $\sigma: \Pi \to \{\uparrow,\to\}$ there
corresponds a cone $C(\sigma)$, and $Pol=\cup_{\sigma} C(\sigma)$.
The intersection of all these cones consists of functions which
are affine on every elementary octahedron (as in Example
3).\medskip

For a function  $f$ on $\Pi$, we denote by $f_1=res_1(f)$ the
restriction of  $f$ to the points of $\Pi$ of the form $(*,*,1)$
(that is on the level 1). Specifically,  $f_1$ is defined on the
grid $\{1,3,\ldots, 2n-1\}\times \{1, 3,\ldots ,2n-1\}$. In
contrast to the algebraic case, for a $T$-polarized function $f$,
the function $f_1$ is not arbitrary, it is supermodular. Let us
recall that a function $g:\mathbb Z^2\to\mathbb R\cup\{-\infty\}$
is called {\em supermodular} if
\begin{equation}\label{superm}
g(i,j)-g(i-1,j)-g(i,j-1)+g(i-1,j-1)\ge 0
\end{equation}
for all $i$ and $j$.
\medskip

{\bf Proposition 1.} {\em If $f$ is a $T$-polarized function then
$f_1$ is supermodular.}\medskip

Indeed, by the definition of a polarized function, we have
                                               $$
f(i-1,j-1,1)+f(i+1,j+1,1)\ge f(i-1,j+1,1)+f(i+1,j-1,1)). \ \ \Box
$$
(It is clear that the restriction of $f$ to each level of $\Pi$ is
a supermodular function as well.) Denote by $Supmod$ the set of
supermodular functions on the grid $\{1,3,\ldots, 2n-1\}\times
\{1, 3,\ldots ,2n-1\}$.\medskip

    {\bf  Proposition 2.} {\em  The mapping $res_1:Pol \to Supmod$
is surjective.  }\medskip

{\em Proof}. Let  $b$ be a supermodular function on
$\{1,3,...,2n-1\}\times \{1,3,...,2n-1\}$. Define the function $f$
on $\Pi$ by the rule:
$$
    f(i,j,k)=b(i-k+1,j-k+1)+b(i-k+3,j-k+3)+...+b(i+k-1,j+k-1).
$$
For example, for $k=1$, this sum consists of a single summand
$b(i,j)$, so that $f_1 =b$. It remains to verify that  $f$ is
$T$-polarized. We do this by proving two claims.\medskip

{\em Claim 1.}
           $$
      f(i-1,j-1,k)+f(i+1,j+1,k)\ge f(i-1,j+1,k)+f(i+1,j-1,k).
             $$
The left hand side is (by definition) the following sum
$b(i-1-k+1,j-1-k+1)+...+b(i-1+k-1,j-1+k-1)+
b(i+1-k+1,j+1-k+1)+...+b(i+1+k-1,j+1+k-1).$ Since $b$ is
supermodular, we have
$$b(i-k,j-k)+b(i-k+2,j-k+2)\ge b(i-k,j-k+2)+b(i-k+2,j-k),$$
$$.......$$
$$b(i+k-2,j+k-2)+b(i+k,j+k)\ge b(i+k-2,j+k)+b(i+k,j+k-2).$$
Summing up these inequalities, on the left hand side we get the
above sum, and on the right hand side we get the sum
$$b(i-k,j-k+2)+...+b(i+k-2,j+k)+$$
$$ b(i-k+2,j-k)+...+b(i+k,j+k-2).$$
The first of these sums is equal to $f(i-1,j+1,k)$ and the second
one is equal to  $f(i+1,j-1,k)$. So, this claim is proven.\medskip

{\em Claim 2.}
                                                            $$
f(i-1,j-1,k)+f(i+1,j+1,k)=f(i,j,k-1)+f(i,j,k+1).
                                                              $$
In fact, both sides are equal to
$b(i-k,j-k)+2b(i-k+2,j-k+2)+...+2b(i+k-2,j+k-2)+b(i+k,j+k)$.
$\Box$\medskip

We see from the proof that the mapping $res_1$ provides a (linear)
bijection between the cone $C(\uparrow,\uparrow,\ldots, \uparrow)$
and the cone $Supmod$. For other cones in $Pol$ this is not the
case. Thus, in general, the `matrix' $f_1$ does not determine
$T$-polarized function $f$.

In the next section we show that a tropical variant of the
`ontogenesis' proves to be more relevant.

\section{Tropical ontogenesis}

We consider a tropical version of the construction of functions
$F$ via genetic arrays $W$ (see Section 2). Let $s=(s(i,j)$, $1\le
i,j\le n$) be a (genetic) array of real numbers (we allow negative
entries as well). We associate to $s$ an (ontogenetic) function
$f=\Phi (s)$ on $\Pi$ by the rule
$$
f(i,j,k)=\max_{(\gamma _1,...,\gamma_k)}[s(\gamma _1)+...+s(\gamma
_k)],
$$
where $(\gamma _1,...,\gamma _k)$ runs over the set of  $k$-tuples
of non-intersecting paths from  $(i-k)/2+1$ to $(j-k)/2+1$,
$\ldots$, from $(i-k)/2+k$ to $(j-k)/2+k$.\medskip

For example,
$f(2,2,1)=\max\{s(1,2)+s(1,1)+s(2,1),s(1,2)+s(2,2)+s(2,1)\}$.\medskip

     {\bf Example  1$'$.} Consider the diagonal array
$s(i,j)=\delta _{ij}$. Then  $\Phi (\delta )=q$, where $q$ is the
function form Example 1. \medskip

     {\bf Example 2$'$.} Consider the array $s\equiv 1$.
Then $\Phi (s)=p$, where $p$ is the function from Example 2. In
fact, any path from  $i$-th column  to $j$-th row contains $i+j-1$
nodes. Therefore, for such a path $\gamma$, $s(\gamma)=i+j-1$.
\medskip

{\bf Theorem 2.} {\em The function  $f=\Phi (s)$ is
$T$-polarized.}\medskip

{\em Proof.} It is a tropicalization of proof of Theorem 1. $\Box$
\medskip

To demonstrate a machinery behind this theorem we consider in
detail one simple particular case.\medskip

     {\bf Example 4.} Consider the elementary octahedron centered
at the point $(1,1,2)$. We would like to check the corresponding
octahedron equality
                                    $$
   f(2,2,2)+f(4,4,2)=\max(f(3,3,1)+f(3,3,3), f(2,4,2)+f(4,2,2)).
                      $$
To start with, we write the values of $f$ as sums of relevant
$s(i,j)$. For the sake of brevity, we will write  $ij$ instead of
$s(i,j)$, and we set $S= 11+12+13+21+22+23+31+32+33 . $\medskip

$f(2,2,2)=11+12+21+22$,

$f(4,4,2)=\max(S-33,S-22,S-11),$

$f(3,3,1)=\max(12+11+21,12+22+21),$

$f(3,3,3)=S,$

$f(2,4,2)=S-31-32-33,$

$f(4,2,2)=S-13-23-33.$\medskip

Thus we have to verify the equality
$$
11+12+21+22+\max(S-33,S-22,S-11)=
$$
$$
\max(\max(12+11+21,12+22+21)+S, S-31-32-33+S-13-23-33).
$$
Subtracting $S$, we are reduced to verifying the relation
$$11+12+21+22+\max(-33,-22,-11)=$$
$$\max(\max(12+11+21,12+22+21), S-31-32-33-13-23-33).$$
Subtracting $11+12+21+22$, we arrive at the obvious equality
$\max(-33,-22,-11)=\max(\max(-22,-11),-33)$. $\Box$\medskip

Denote by $Arr=\mathbb{R}\otimes (\{1,...,n\}\times
\{1,...,n\})\cong\mathbb{R}^{n\times n}$ the set of  $n\times n$
arrays. Due to Theorem 2, we have a piecewise linear mapping
$$
                         \Phi :Arr \to Pol.
$$

{\bf Proposition 3.} {\em The mapping $ \Phi$  is a bijection.}
\medskip

{\em Proof}. We construct the inverse mapping. Let $f$ be a
$T$-polarized function, and let $g=\alpha^*(f)$. Since $g$
determines $f$ (see Lemma 1 below), the proposition will follow
from a bijection between $s$ and $g$. By definition of $\Phi $,
$g(i,j)$ is equal to $\sum_{a\le i,b\le j} s(a,b)$. Therefore, in
turn,
$$
           s(i,j)=g(i,j)-g(i-1,j)-g(i,i-1)+g(i-1,j-1). \ \ \Box
$$

Note that  $s$ and $g$ are related by an invertible linear
transformation. \medskip

Thus we see that any $T$-polarized function  $f$ is determined by
its genetic array  $s$. As we have seen, the genotype $s\equiv 1$
determines the function $p$ from Example 2.

One can consider a more general case. Let $MD$ be the set of
arrays $s$ such that $s(i,j)\le s(i+1,j+1)$ for all $i$, $j$ (when
both sides are defined). Then, for any $(i,j,k)\in \Pi$, we have
\[
\Phi(s)(i,j,k)=\sum_{2a\le i+k,\, 2b\le j+k}s(a,b)- \sum_{2a'\le
i-k,\, 2b'\le j-k}s(a',b').
\]
Indeed, for an array $s\in MD$, one can take the maximal path from
column $i$ to row $j$ in the shape of a convex hook.

\section{$RSK$-correspondence}

In Section 2 we considered two maps  $\alpha$ and  $\beta$ from
the square grid $\{0,1,...,n\}\times \{0,1,...,n\}$ to  the
pyramid $\Pi$. As in Section 2, one can use these maps to restrict
$T$-polarized functions on this grid. Note that  the restricted
functions vanish on the south-west boundary of this grid. So, we
actually  deal with the functions on the grid of size $n\times
n$.\medskip

{\bf Lemma 1}. {\em The maps} $\alpha^*  : Pol \to \mathbb{R}
^{n\times n}$ and $\beta^*  : Pol \to \mathbb{R} ^{n\times n}$
{\em are bijections}.\medskip

{\em Proof}. This follows from the octahedron recurrence
(\ref{trop}). Let us check that  $\alpha^*$ is a surjection. Using
(\ref{trop}), we propagate a function $g$ from faces $ABE$ and
$ACE$) to points of $\Pi$. That is, for every elementary
octahedron such that the function is defined at five of its
vertices except the end point of the diagonal parallel to the
propagation vector $(2,2,0)$, we set value at the remaining vertex
using (\ref{trop}). Suppose that as a result of this procedure
some of the points of $\Pi$ remain non-filled. Let  $p$ be such a
point with minimal value of $i+j$. This point can not lie in the
base of the pyramid or in the faces $ABE$ and $ACE$ since all
these are valued. Therefore, we can move from $p$ in $\Pi$ along
the vectors $(-2,-2,0)$, $(0,-2,0)$, $(-2,0,0)$, $(-1,-1,1)$, and
$(-1,-1,1)$. On the other hand, the function $f$ is already
defined at these nodes; hence, by (\ref{trop}), we can define
$f(p)$.

The injectivity of $\alpha^*$ is proven similarly.

The case of the map  $\beta^*$  is dealt with in a similar way.
 $\Box$\medskip

Thus we have three bijections
$$
\begin{CD}
\mathbb{R}^{n\times n} @<{\alpha^*}<< Pol @>{\beta^*}>>
\mathbb{R}^{n\times n}
\\
@. @A{\Phi}AA @. \\
{} @. Arr @. {}
\end{CD}
$$

Now we characterize the sets  $\alpha^*(\Phi(Arr_+))$ and
$\beta^*(\Phi(Arr_+))$, where $Arr_+\cong \mathbb R_+^{n\times n}$
is the set of non-negative arrays. The OR-map $\beta^* \circ
(\alpha ^*)^{-1}$ gives a bijection between these sets.

The image $\alpha^*(\Phi(Arr_+))$ admits a rather simple
description. Namely, it consists of supermodular functions which
vanish on the south-west boundary of the $(n+1)\times (n+1)$-grid.
Indeed, if $g=\alpha^*(\Phi(s))$ then
$s(i,j)=g(i,j)+g(i-1,j-1)-g(i,j-1)+ g(i-1,j)$ by Proposition 3.
Therefore $g$ is supermodular if and only if $s$ is non-negative.

To describe $\beta^*(\Phi(Arr_+))$, we need a notion complementary
to that of supermodularity. A function $h$ on a square (or
rectangular) grid is called {\em inframodular} if
$$
               h(i,j)+h(i+1,j)\ge h(i,j-1)+h(i+1,j+1)
$$
and
$$
               h(i,j)+h(i,j+1)\ge h(i-1,j)+h(i+1,j+1)
$$
for all $i$ and $j$ (when all the terms are defined). The
inframodularity means that the ``diagonal partial difference''
$\partial _dh(i,j):=h(i,j)-h(i-1,j-1)$ is decreasing as a function
of $(i,j)$.

A function is called {\em discretely concave} if it is
supermodular and inframodular.

We claim  that the image $\beta^*(\Phi(Arr_+))$ consists of
inframodular functions. More precisely, we have the
following\medskip

     {\bf Theorem 3.} {\em Let $g$ be a supermodular function on
$\{0,1,...,n\}\times \{0,1,...,n\}$ such that $g(*,0)=g(0,*)=0$.
Then the function $h=\beta^*(\alpha^{*-1}(g))$ is inframodular,
$h(*,0)=h(0,*)=0$, and $h(n-1,n-1)\le h(n,n)$. The converse is
also true.}\medskip

We prove Theorem 3 in the next section.\medskip

While supermodular functions are related to non-negative arrays,
inframodular functions are related to plane partitions. Here by
{\em plane partition} we mean an arbitrary weakly decreasing
function $p: \{1,...,n\}\times \{1,...,n\} \to \mathbb R_+$. If we
replace $\mathbb{R}_+$ on $\mathbb{Z}_+$, we obtain classical
plane partitions.

Now, if we have an inframodular function $h$ on
$\{0,1,...,n\}\times \{0,1,...,n\}$ then $p=\partial_d (h)$ is a
weakly decreasing function on $\{1,...,n\}\times \{1,...,n\}$. If,
in addition, $h(n,n)\ge h(n-1,n-1)$ then $p(n,n)\ge 0$ from which
it follows that the same holds for all $p(i,j)$. Finally, the
function $p$ uniquely determines the boundary function $h$
provided that $h$ vanishes on the south-west boundary. To wit,
$h(i,j)=p(i,j)+p(i-1,j-1)+p(i-2,j-2)+...$ .\medskip

Thus, due to Theorem 3, we get the following commutative diagram
of bijections\medskip

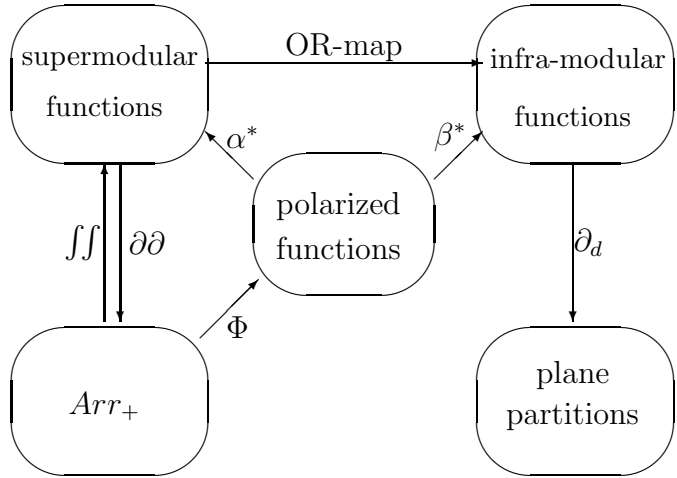
\unitlength=0.7mm \special{em:linewidth 0.4pt}
\linethickness{0.4pt}
\begin{picture}(143.00,110.00)(-10,20)
\put(38.00,45.00){\oval(37.00,28.00)[]}
\put(38.00,105.00){\oval(37.00,30.00)[]}
\put(82.00,78.50){\oval(34.00,27.00)[]}
\put(125.50,105.00){\oval(37.00,30.00)[]}
\put(125.50,45.00){\oval(37.00,28.00)[]}
\put(65.00,87.00){\vector(-1,1){9.00}}
\put(99.00,87.00){\vector(1,1){9.00}}
\put(56.00,109.00){\vector(1,0){52.00}}
\put(37.00,60.00){\vector(0,1){30.00}}
\put(125.00,90.00){\vector(0,-1){30.00}}
\put(40.00,90.00){\vector(0,-1){30.00}}

\put(55.00,57.00){\vector(1,1){11.00}}
\put(37.00,44.00){\makebox(0,0)[cc]{$Arr_+$}}
\put(62.00,59.00){\makebox(0,0)[cc]{$\Phi$}}
\put(81.00,82.00){\makebox(0,0)[cc]{polarized}}
\put(81.00,74.00){\makebox(0,0)[cc]{functions}}
\put(38.00,110.00){\makebox(0,0)[cc]{\small supermodular}}
\put(37.00,101.00){\makebox(0,0)[cc]{\small functions}}
\put(33,75.00){\makebox(0,0)[cc]{$\int\!\!\int$}}
\put(45,75.00){\makebox(0,0)[cc]{$\partial\partial$}}
\put(63.00,95.00){\makebox(0,0)[cc]{$\alpha^*$}}
\put(102.00,95.00){\makebox(0,0)[cc]{$\beta^*$}}
\put(82.00,112.00){\makebox(0,0)[cc]{OR-map}}
\put(126.00,110.00){\makebox(0,0)[cc]{\small infra-modular}}
\put(125.00,99.00){\makebox(0,0)[cc]{\small functions}}
\put(128.00,75.00){\makebox(0,0)[cc]{$\partial_d$}}
\put(125.00,50.00){\makebox(0,0)[cc]{plane}}
\put(125.00,42.00){\makebox(0,0)[cc]{partitions}}
\end{picture}

\noindent In particular, the composition $\partial_d \circ \beta^*
\circ \Phi$ gives a natural bijection between the set of
non-negative arrays and that of plane partitions which can be
considered as a kind of RSK-correspondence (see \cite{umn},
14.6-14.8). Since octahedron recurrence (\ref{trop}) propagates
integer-valued data to integer-valued data, our bijection  sends
non-negative integer arrays to ordinary integer-valued plane
partitions.\medskip

We conclude this section by a comparison of the bijection
$\partial_d \circ \beta^* \circ \Phi$ with the classical
RSK-correspondence. This relationship is based on a simple natural
bijection between the set of plane partitions and the set of pairs
of semi-standard Young tableaus (with $n$ rows) of the same shape.
Let us recall the definition of this bijection (see the details in
\cite{Stan-2} or \cite{umn}).

Let $p$ be an (integer-valued) plane partition. Then the following
weakly decreasing $n$-tuple $\lambda=(p(1,1),p(2,2),...,p(n,n))$
is a partition (of the number $f(n,n,n)=\sum_{a,b} s(a,b)$ when
$p$ originates from genetic array $s$. This partition $\lambda$ is
the shape of the array $s$ in the sense of \cite{umn}).

Now consider the $(n-1)$-tuple
$\lambda'=(p(2,1),p(3,2),...,p(n,n-1))$ ($\lambda'$ is the shape
of the array $s'$ obtained from $s$ by forgetting the last
column). Because of monotonicity, the sequences $\lambda$ and
$\lambda '$ interlace, that is
$$
                   \lambda_1 \ge \lambda'_1 \ge \lambda_2 \ge ...\ge
\lambda_{n-1}\ge \lambda'_{n-1}\ge \lambda_n .
                  $$
The same holds for the other diagonals. Thus, the restriction of
$p$ on the triangle below the principal diagonal yields a sequence
of $n$ interlacing partitions $(\lambda^{(0)}=\lambda,
\lambda^{(1)}=\lambda', \lambda^{(2)},...,
\lambda^{(n)}=\emptyset)$, that is a semi-standard Young tableau.
Specifically, to obtain the tableau one needs to fill every
horizontal strip $\lambda^{(i)}-\lambda^{(i+1)}$ with the letter
$n-i$.

Analogously, the upper half of $p$ gives another Young tableau.
This tableau has the same shape $\lambda$. Conversely, a pair of
tableaux of the same shape gives a plane partition.

This permits to reformulate Theorem 3 as follows\medskip

{\bf Theorem 3$'$. } {\em The map $\partial_d \circ \beta ^* \circ
\Phi $ gives a bijection between non-negative integer $(n\times
n)$-arrays and pairs of semi-standard  Young tableaux of the same
shape (filled with the alphabet $\{1,\ldots ,n\}$). }\medskip

It is easy to see (using Appendix B from \cite{umn}) that the
bijection coincides with the (modified) RSK-correspondence from
\cite{umn}. To order to get the original RSK, we should replace
one of the tableaux, viz `$Q$-symbol', which corresponds to the
half of the function $h$ living on the face $BED$, by its
Sch\"utzenberger transform.\medskip

{\bf Example 5.} We illustrate this theorem by an example.
Consider an array $s={\small
\begin{array}{ccc}
  1 & 2 & 3 \\
  1 & 1 & 5 \\
  2 & 3 & 1
\end{array}.}$
The mapping $\int\!\!\int$ takes $s$ to the function $g={\small
\begin{array}{cccc}
 0 & 4 & 10 & 18 \\
 0 & 3 & 7 & 13 \\
 0 & 2 & 5 & 6 \\
 0 & 0 & 0 & 0
\end{array}}$. Using the octahedron recurrence, we propagate $g$ to
$\Pi$ (see the picture below)

\unitlength=.5mm \special{em:linewidth 0.4pt}
\linethickness{0.4pt}
\begin{picture}(109.00,140.00)(-50,-10)
\bezier{40}(15.00,20.00)(37.00,65.00)(60.00,110.00)
\bezier{24}(15.00,20.00)(45.00,13.00)(75.00,5.00)
\bezier{13}(75.00,5.00)(90.00,13.00)(105.00,20.00)
\bezier{42}(75.00,5.00)(67.00,55.00)(60.00,110.00)
\bezier{40}(60.00,110.00)(82.00,65.00)(105.00,20.00)
\bezier{24}(45.00,35.00)(72.00,28.00)(105.00,20.00)
\bezier{13}(45.00,35.00)(30.00,28.00)(15.00,20.00)
\put(30.00,50.00){\makebox(0,0)[cc]{2}}
\put(30.00,50.00){\makebox(0,0)[cc]{2}}
\put(50.00,45.00){\makebox(0,0)[cc]{5}}
\put(70.00,40.00){\makebox(0,0)[cc]{6}}
\put(40.00,55.00){\makebox(0,0)[cc]{3}}
\put(50.00,60.00){\makebox(0,0)[cc]{4}}
\put(45.00,80.00){\makebox(0,0)[cc]{7}}
\put(55.00,85.00){\makebox(0,0)[cc]{10}}
\put(65.00,75.00){\makebox(0,0)[cc]{13}}
\put(60.00,50.00){\makebox(0,0)[cc]{6}}
\put(70.00,55.00){\makebox(0,0)[cc]{7}}
\put(80.00,45.00){\makebox(0,0)[cc]{8}}
\put(90.00,50.00){\makebox(0,0)[cc]{11}}
\put(75.00,80.00){\makebox(0,0)[cc]{17}}
\put(60.00,110.00){\makebox(0,0)[cc]{18}}
\put(15.00,20.00){\makebox(0,0)[cc]{0}}
\put(25.00,25.00){\makebox(0,0)[cc]{0}}
\put(35.00,30.00){\makebox(0,0)[cc]{0}}
\put(45.00,35.00){\makebox(0,0)[cc]{0}}
\put(35.00,15.00){\makebox(0,0)[cc]{0}}
\put(45.00,20.00){\makebox(0,0)[cc]{0}}
\put(55.00,25.00){\makebox(0,0)[cc]{0}}
\put(65.00,30.00){\makebox(0,0)[cc]{0}}
\put(55.00,10.00){\makebox(0,0)[cc]{0}}
\put(65.00,15.00){\makebox(0,0)[cc]{0}}
\put(75.00,20.00){\makebox(0,0)[cc]{0}}
\put(85.00,25.00){\makebox(0,0)[cc]{0}}
\put(75.00,5.00){\makebox(0,0)[cc]{0}}
\put(85.00,10.00){\makebox(0,0)[cc]{0}}
\put(95.00,15.00){\makebox(0,0)[cc]{0}}
\put(105.00,20.00){\makebox(0,0)[cc]{0}}
\put(9.00,17.00){\makebox(0,0)[cc]{$A$}}
\put(75.00,-2.00){\makebox(0,0)[cc]{$B$}}

\put(111.00,18.00){\makebox(0,0)[cc]{$D$}}
\put(42.00,40.00){\makebox(0,0)[cc]{$C$}}
\put(66.00,117.00){\makebox(0,0)[cc]{$E$}}
\bezier{30}(60.00,110.00)(53.00,72.00)(45.00,35.00)
\bezier{16}(30.00,50.00)(50.00,45.00)(70.00,40.00)
\bezier{8}(70.00,40.00)(80.00,45.00)(90.00,50.00)
\bezier{16}(90.00,50.00)(70.00,55.00)(50.00,60.00)
\put(50.00,60.00){\makebox(0,0)[cc]{4}}
\put(55.00,85.00){\makebox(0,0)[cc]{10}}
\bezier{8}(30.00,50.00)(40.00,55.00)(50.00,60.00)
\bezier{8}(45.00,80.00)(55.00,78.00)(65.00,75.00)
\bezier{4}(65.00,75.00)(70.00,77.00)(75.00,80.00)
\bezier{8}(75.00,80.00)(65.00,83.00)(55.00,85.00)
\bezier{4}(55.00,85.00)(50.00,83.00)(45.00,80.00)
\end{picture}

\

\noindent On the faces $BED\cup CED$ we get the function
$h={\small
\begin{array}{cccc}
  0 & 6 & 13 & 18 \\
  0 & 8 & 17 & 10 \\
  0 & 11 & 7 & 4 \\
  0 & 0 & 0 & 0
\end{array}.}$
Applying $\partial_d$ gives the plane partition $p={\small
\begin{array}{ccc}
   6 & 5 & 1 \\
   8 & 6 & 3 \\
   11 & 7 & 4
\end{array}.}$
 The lower half of $p$ gives the Young tableau\medskip

3

2\ 2\ 2\ 3\ 3\ 3

1\ 1\ 1\ 1\ 2\ 2\ 2\ 3\ 3\ 3\ 3\medskip

\noindent of the shape $(11,6,1)$. The upper half of $p$ gives the
Young tableau\medskip

3

2\ 2\ 2\ 2\ 2\ 3

1\ 1\ 1\ 1\ 1\ 1\ 2\ 2\ 3\ 3\ 3\medskip

\noindent of the same shape $(11,6,1)$.\medskip

\section{Proof of Theorem 3}

Before proving Theorem 3, we would like to illustrate it by a
simple example. Consider the function $g$ given by
$$\begin{array}{ccc}
  0 & c & e \\
  0 & a & b \\
  0 & 0 & 0
\end{array}.
$$
Then its OR-image $h$ is
$$\begin{array}{ccc}
  0 & b & e \\
  0 & a' & c \\
  0 & 0 & 0
\end{array},$$
where (according to (\ref{trop})) $a+a'=\max(e+0,b+c)$.

If $g$ is supermodular, then the following four inequalities hold
$$
                        b,c\ge a\ge 0 \text{ and } a+e\ge b+c.
$$
In order to verify the inframodularity of $h$, we have to check
the following five inequalities
$$
                  a'\ge b,c; \ a'+b\ge e;\  a'+c\ge e; \  e\ge a'.
$$
Since $b+c\ge b+a$, we have $a'+a=\max(e,b+c)\ge b+a$ and,
consequently, $a'\ge b$. Since $b\ge a$, we have
$a'+a+b=\max(e,b+c)+b\ge e+b\ge e+a$; cancellation of $a$ gives
$a'+b\ge e$. Finally, since $e+a\ge e$ and $e+a\ge b+c$, we have
$e+a\ge \max(e+0,b+c)=a'+a$ and $e\ge a'$.

Vice versa, let  $h$ satisfy these five inequalities. Then adding
$a$ to both sides of the inequality $e\ge a'$, we get  $e+a\ge
\max(e,b+c)$. This implies $e+a\ge e$, and hence  $a\ge 0$ and
$e+a\ge b+c$. Since $a'+b\ge e$ and $a'+b\ge c+b$, we have
$a'+b\ge \max(e,b+c)=a'+a$, that is  $b\ge a$. Thus $g$ is
supermodular.\medskip

{\em Proof of Theorem 3}.  We will exploit a result from
\cite{arxiv}. For this we need to transform our data in order to
satisfy the conditions of (\cite{arxiv}, Corollary 5). To wit, we
have to transform data from the pyramid to a prism.

Let us turn the face  $ACE$ around the edge  $AE$ in order to get
a prism $AC'CBED$. We locate the function $g$ on the front face
$AC'EB$, and we assign zero values to the points of the triangle
face $AC'C$. In accordance with (\ref{trop}), we propagate these
data to the prism. In view of the non-negativity of $g$, the
values of the propagated function on $ACE$ coincide with the old
values. Therefore the octahedron recurrence on the prism gives the
same polarized function  $f$ on the pyramid $ABCDE$.

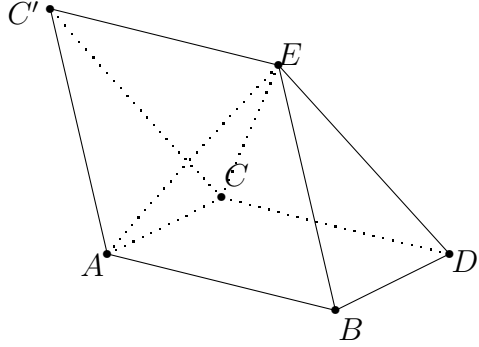
\unitlength=.5mm \special{em:linewidth 0.4pt}
\linethickness{0.4pt}
\begin{picture}(134.00,100.00)(-45,-5)
\put(40.00,20.00){\circle*{2.00}}
\put(100.00,5.00){\circle*{2.00}}
\put(130.00,20.00){\circle*{2.00}}
\put(70.00,35.00){\circle*{2.00}}
\put(85.00,70.00){\circle*{2.00}}
\put(25.00,85.00){\circle*{2.00}}
\emline{25.00}{85.00}{1}{40.00}{20.00}{2}
\emline{40.00}{20.00}{3}{100.00}{5.00}{4}
\emline{100.00}{5.00}{5}{85.00}{70.00}{6}
\emline{85.00}{70.00}{7}{130.00}{20.00}{8}
\emline{130.00}{20.00}{9}{100.00}{5.00}{10}
\emline{85.00}{70.00}{11}{25.00}{85.00}{12}
\bezier{13}(40.00,20.00)(55.00,27.00)(70.00,35.00)
\bezier{24}(70.00,35.00)(100.00,27.00)(130.00,20.00)
\bezier{26}(70.00,35.00)(50.00,56.00)(25.00,85.00)
\put(36.00,17.00){\makebox(0,0)[cc]{$A$}}
\put(104.00,0){\makebox(0,0)[cc]{$B$}}
\put(74.00,41.00){\makebox(0,0)[cc]{$C$}}
\put(134.00,18.00){\makebox(0,0)[cc]{$D$}}
\put(88.00,73.00){\makebox(0,0)[cc]{$E$}}
\put(18.00,84.00){\makebox(0,0)[cc]{$C'$}}
\bezier{26}(40.00,20.00)(61.00,45.00)(85.00,70.00)
\bezier{15}(70.00,35.00)(77.00,52.00)(85.00,70.00)
\end{picture}

Now, we can improve negative $H$-breaks, $ g(i,j)+g(i+1,j)-
g(i,j-1)-g(i+1,j+1)$ and negative $V$-breaks, $g(i,j)+g(i,j+1)-
g(i-1,j)-g(i+1,j+1)$, of the function $g$ by adding to $g$ two
appropriate functions of the variables $i+k$ and $j+k$
respectively. As a result we get that the modified function
$\tilde g$ on $ABEC'$ (as well as the corresponding function on
$AC'C$) is a discrete concave function (that is supermodular and
inframodular). By Corollary 5 of Theorem 1 in \cite{arxiv}, the
corresponding polarized function $\tilde f$ is a discrete concave
function on $\Pi$. In particular, $\tilde f$ is discretely concave
on the faces $DBE$ and $DCE$. Now we subtract the functions which
we added in order to improve negative $H$-, $V$-breaks of $g$. We
obtain the `old' function $h$. One can see that this subtraction
does not change the $H$- and $V$-breaks on the boundary $BED\cup
CED$. Thus, the function $h$ is infra-modular.

It remains to verify $h(n-1,n-1)\le h(n,n)$, that is
$a'=f(n+1,n+1,n-1)\le f(n,n,n)=e$. By induction we can suppose
that $e'=f(n,n,n-2)\le f(n-1,n-1,n-1)=a$. Set  $b=f(n+1,n-1,n-1)$
and $c=f(n-1,n+1,n-1)$. Because of the octahedron recurrence we
get
$$
a'+a=\max(e+e', b+c).
$$
Supermodularity of  $g$ implies  $a+e\ge b+c$. Therefore
$$
          a'+a\le \max(e+e',a+e)=\max(e',a)+e\le \max(a,a)+e=a+e,
$$
and we get  $a'\le e$.

The converse is proved by a similar construction. First we turn
the face  $BED$ around the edge $DE$ and get a square $CEB'D$. We
locate the function $h$ on this square. We set it to be
identically zero on the new face  $BB'D$. These data are
propagated by the octahedron recurrence (\ref{trop}) in the
opposite direction $(-2,-2,0)$. On the pyramid, we obtain the
function $f$. Again, we can add to $h$ an appropriate function of
$i+k$ in order to improve negative supermodular breaks of $h$. We
get a discrete concave function $\tilde h$. Therefore (by the same
Corollary 5 from \cite{arxiv}) the corresponding $\tilde f$ is a
discrete concave function on the prism. Subtracting from this
function the addendum to $h$, we come back to $f$. This
subtraction does not change supermodular breaks of $g$ except,
possibly, the supermodular breaks along the main diagonal.

Thus, to complete the proof, we have to show that $g$ has
non-negative breaks along the main diagonal. That is (in the above
notations) we need to show that $a+e\ge b+c$. We have $a'\le e$
and
$$
                       a'+a=\max(e+e', b+c).
$$
Hence
$$
               a+e=a+a'+e-a'\ge a+a'=\max(e+e',b+c)\ge b+c.
            $$
Moreover, $a+e\ge e+e'$, which proves the inequality $a\ge e'$ and
yields a basis for an induction argument. This completes the proof
of supermodularity of $g$ and that of the theorem.

\end{document}